\newcommand{\al}{\alpha}               
\newcommand{\be}{\beta}
               
\newcommand{\Ga}{\Gamma}
               
\newcommand{\De}{\Delta}
\newcommand{\lb}{\lambda}              

\newcommand{\sig}{\sigma}

\newcommand{\vphi}{\varphi}

\newcommand{\cal}{\mathcal}

\newcommand{\calf}{{\cal F}}

\newcommand{\calm}{{\cal M}}           

\newcommand{\calr}{{\cal R}}

\newcommand{\calv}{{\cal V}}

\newcommand{\Fix}{{\rm Fix}}

\newcommand{\incl}{\subseteq}         
\newcommand{\aincl}{\supseteq}
\newcommand{\sincl}{\subset}          
\newcommand{\saincl}{\supset}
\newcommand{\es}{\emptyset}           
\newcommand{\sm}{\setminus}
       
\newcommand{\limpl}{\Longrightarrow}
   
\newcommand{\lequi}{\Longleftrightarrow}
\newcommand{\oo}{\infty}

\newcommand{\sk}{\smallskip}
       
\newcommand{\n}{\noindent}

\def\R+oo{R_+\cup\{\oo\}}

\def\dtends   {\stackrel {\it d}{\longrightarrow}}

\def\(V)tends {\stackrel {(\calv)}{\longrightarrow}}

\newcommand{\barr}{\begin{array}}        
\newcommand{\earr}{\end{array}}
\newcommand{\bcor}{\begin{corollary}}    
\newcommand{\ecor}{\end{corollary}}
\newcommand{\ben}{\begin{enumerate}}     
\newcommand{\een}{\end{enumerate}}
\newcommand{\beq}{\begin{equation}}       
\newcommand{\eeq}{\end{equation}}
\newcommand{\bex}{\begin{example}}        
\newcommand{\eex}{\end{example}}
\newcommand{\blemma}{\begin{lemma}}       
\newcommand{\elemma}{\end{lemma}}
\newcommand{\bproof}{\begin{proof}}       
\newcommand{\eproof}{\end{proof}}
\newcommand{\bprop}{\begin{proposition}}  
\newcommand{\eprop}{\end{proposition}}
\newcommand{\brem}{\begin{remark}}        
\newcommand{\erem}{\end{remark}}
\newcommand{\btab}{\begin{tabular}}       
\newcommand{\etab}{\end{tabular}}
\newcommand{\btheorem}{\begin{theorem}}   
\newcommand{\etheorem}{\end{theorem}}

\documentclass[reqno]{amsart}

\newtheorem{theorem}{\bf Theorem}
\newtheorem{corollary}{\bf Corollary}
\newtheorem{example}{\bf Example}
\newtheorem{lemma}{\bf Lemma}
\newtheorem{proposition}{\bf Proposition}
\newtheorem{remark}{\bf Remark}

\begin{document}

\title
[Ultrametric Fixed Points in Reduced Axiomatic Systems]
{ULTRAMETRIC FIXED POINTS IN \\
REDUCED AXIOMATIC SYSTEMS}

\author{Mihai Turinici}
\address{
"A. Myller" Mathematical Seminar;
"A. I. Cuza" University;
700506 Ia\c{s}i, Romania
}
\email{mturi@uaic.ro}




\begin{abstract}
The Brezis-Browder ordering principle
[Advances Math., 21 (1976), 355-364]
is used to get a proof,  
in the reduced axiomatic system
(ZF-AC+DC),
of a fixed point result  
[in the complete axiomatic system (ZF)]
over Cantor complete ultrametric spaces 
due to Petalas and Vidalis
[Proc. Amer. Math. Soc., 118 (1993), 819-821]. 
\end{abstract}

\maketitle

\section{Introduction}
\setcounter{equation}{0}

Throughout this exposition, 
the axiomatic system in use is 
Zermelo-Fraenkel's (in short: ZF), 
as described by 
Cohen \cite[Ch 2]{cohen-1966}.
The notations and basic facts about 
its axioms are more or less usual.

Remember that, an outstanding part of it is 
the 
{\it Axiom of Choice} (abbreviated: AC);
which, in a convenient manner, may be written as
\ben 
\item[] (AC)\ \ 
For each nonempty set $X$, there exists 
a (selective) function \\
$f:(2)^X\to X$
with $f(Y)\in Y$, for each $Y\in (2)^X$.
\een
[Here, $(2)^X$ denotes the class of
all nonempty parts in $X$].
There are many logical equivalents of (AC);
see, for instance, 
Moore \cite[Appendix 2]{moore-1982}.
A basic one is the
{\it Zorn-Bourbaki Maximal Principle}
(in short: ZB), expressed as
\ben
\item[] (ZB)\ \  
Let the partially ordered set $(X,\le)$ be 
{\it inductive}
(any totally ordered part $C$ of $X$ 
is {\it bounded above}:
$C\le b$ (i.e.: $x\le b$, $\forall x\in C$), 
for some $b\in X$).
Then, for each (starting) $u\in X$, 
there exists a maximal element 
$v\in X$ \\
(in the sense: $v\le z\in X$ implies $v=z$),
with $u\le v$;
\een
for a direct proof of this
(avoiding transfinite induction),
see 
Bourbaki \cite{bourbaki-1949}.

Let $X$ be a nonempty set.
By a {\it sequence} in $X$, we mean any mapping $x:N\to X$;
where $N:=\{0,1,...\}$ is the set of {\it natural} numbers.
For simplicity reasons, it will be useful to denote it as 
$(x(n); n\ge 0)$, or $(x_n; n\ge 0)$;
moreover, when no confusion can arise, 
we further simplify this notation as 
$(x(n))$ or $(x_n)$, respectively.
Also, any sequence $(y_n:=x_{i(n)}; n\ge 0)$ with
\ben
\item[]
$(i(n); n\ge 0)$ is {\it divergent}\ 
[i.e.: $i(n)\to \oo$ as $n\to \oo$], 
\een
will be referred to as a
{\it subsequence} of $(x_n; n\ge 0)$.
Call the subset $Y$ of $X$, 
{\it almost singleton} (in short: {\it asingleton})
provided [$y_1,y_2\in Y$ implies $y_1=y_2$];
and {\it singleton}
if, in addition, $Y$ is nonempty;
note that in this case,
$Y=\{y\}$, for some $y\in X$. 
Further, let $d:X\times X\to R_+:=[0,\oo[$ 
be a {\it metric} over $X$;
the couple $(X,d)$ 
will be termed a {\it metric space}.
Finally, let $T\in \calf(X)$ be a selfmap of $X$.
[Here, for each couple $A, B$ of nonempty sets,
$\calf(A,B)$ stands for the class of all functions 
from $A$ to $B$; 
when $A=B$, we write $\calf(A)$ in place of $\calf(A,A)$].
Denote $\Fix(T)=\{x\in X; x=Tx\}$;
each point of this set is referred to as 
{\it fixed} under $T$.
In the metrical fixed point theory, 
such points are to be determined 
according to the context below,
comparable with the one 
described in 
Rus \cite[Ch 2, Sect 2.2]{rus-2001}:

{\bf pic-1)} 
We say that $T$ is a {\it Picard operator} (modulo $d$) if,
for each $x\in X$, 
the iterative sequence
$(T^nx; n\ge 0)$ is $d$-convergent

{\bf pic-2)} 
We say that $T$ is a {\it strong Picard operator} (modulo $d$) if,
for each $x\in X$, $(T^nx; n\ge 0)$ is $d$-convergent 
with $\lim_n(T^nx)\in \Fix(T)$

{\bf pic-3)} 
We say that $T$ is {\it fix-asingleton} 
(resp., {\it fix-singleton}) if
$\Fix(T)$ is asingleton (resp., singleton).

In this perspective, 
a basic answer to the posed question
is the 1922 one, due to
Banach \cite{banach-1922}.
Given $\al\ge 0$, 
let us say that $T$ 
is {\it $(d;\al)$-contractive}, provided 
\ben
\item[] (a01)\ \ 
$d(Tx,Ty)\le \al d(x,y)$,\ for all  $x,y\in X$.
\een

\btheorem \label{t1}
Suppose that  
$T$ is $(d;\al)$-contractive, for some $\al\in [0,1[$.
In addition, let $(X,d)$ be complete.
Then, $T$ is a strong Picard operator 
(modulo $d$), and 
fix-asingleton (hence, fix-singleton). 
\etheorem

This result -- 
referred to as {\it Banach's contraction principle} --
found a multitude of applications in 
operator equations theory;
so, it was the subject of many extensions.
A natural way of doing this  
is by considering "functional" contractive conditions
\ben
\item[] (a02)\ \ 
$d(Tx,Ty)\le F(d(x,y),d(x,Tx),d(y,Ty),d(x,Ty),d(y,Tx))$,\\
for all  $x,y\in X$;
\een
where $F:R_+^5\to R_+$ is a function.
Some important 
results in the area 
have been established by 
Boyd and Wong \cite{boyd-wong-1969},
Matkowski \cite{matkowski-1975},
and
Leader \cite{leader-1979}.
For more details about other possible choices of $F$
we refer to the 1977 paper by
Rhoades \cite{rhoades-1977};
some extensions of these to quasi-ordered structures
may be found in 
Turinici \cite{turinici-1986-DM}.
Further, a natural extension of 
the contractive condition above is
\ben
\item[] (a03)\ \ 
($T$ is {\it $d$-strictly-nonexpansive}):\\
$d(Tx,Ty)< d(x,y)$, for all $x,y\in X$, $x\ne y$.
\een
Note that, a fixed point for such maps is to
be reached when 
$(X,d)$ is a 
{\it compact metric space};
cf. Edelstein \cite{edelstein-1961}.
Another circumstance when this conclusion holds
is that of $(X,d)$ being a
{\it (transfinite) Cantor complete} 
ultrametric space;
see, for instance,
Petalas and Vidalis \cite{petalas-vidalis-1993}.
In this last case, 
a basic tool used in authors' proof is 
(ZB) 
(=the {\it Zorn-Bourbaki Maximal Principle}); 
or, equivalently (see above): 
(AC) 
(=the {\it Axiom of Choice});
hence, this fixed point result is 
valid in the {\it complete} 
Zermelo-Fraenkel system (ZF).
However, since all arguments used there
are countable in nature,
it is highly expectable that 
a denumerable version of (ZB) 
should suffice for the result's conclusion 
to hold.
It is our aim in the present exposition to 
prove that this is indeed the case.
Precisely, we show that the
Zorn-Bourbaki Maximal Principle
appearing there may be replaced with 
a countable version of it --
namely, the Brezis-Browder ordering principle 
\cite{brezis-browder-1976} --
to solve the posed fixed point question;
hence, the Petalas-Vidalis result
is ultimately deductible in
the {\it reduced} Zermelo-Fraenkel system
(ZF-AC+DC); where 
(DC) is the {\it Principle of Dependent Choices}.
Note that, the proposed reasoning is applicable 
as well to many other statements of this type;
such as 
the ones due to 
Mishra and Pant \cite{mishra-pant-2014}.
Further aspects will be delineated elsewhere.

\section{Brezis-Browder principles}
\setcounter{equation}{0}

Let $M$ be a nonempty set. Take a {\it quasi-order} $(\le)$
[i.e.: a
{\it reflexive} ($x\le x$, $\forall x\in X$) 
and 
{\it transitive} ($x\le y$, $y\le z$ $\limpl$ $x\le z$) 
relation] over it;
the pair $(M,\le)$ will be then referred to as a
{\it quasi-ordered structure}.
Let also $\vphi:M\to R_+$ be a function.
Call the point $z\in M$, $(\le,\vphi)$-{\it maximal} when:
$z\le w\in M$  implies $\vphi(z)=\vphi(w)$.
A basic result about  such points is the 1976
Brezis-Browder ordering principle
\cite{brezis-browder-1976} (in short: BB).

\bprop \label{p1}
Suppose that 
the quasi-ordered structure $(M,\le)$ 
and the function $\vphi$ (taken as before) fulfill
\ben
\item[] (b01)\ \   
$(M,\le)$ is sequentially inductive: \\
each ascending sequence has an upper bound (modulo $(\le)$)
\item[] (b02)\ \ 
$\vphi$ is $(\le)$-decreasing ($x\le y \limpl \vphi(x) \ge \vphi(y)$).
\een
Then, for each $u\in M$ there exists a $(\le, \vphi)$-maximal
$v\in M$ with $u\le v$.
\eprop

{\bf (A)}
In particular, assume that (in addition)
\ben
\item[]
$(\le)$ is {\it antisymmetric} 
($x\le y$, $y\le x$ $\limpl$ $x=y$).
\een
We then say that it is a {\it (partial) order} on $M$;
and the pair $(M,\le)$ will be called a
{\it (partially) ordered structure}.
In this case, by an appropriate choice of 
our structure (related to existence of
functions $\vphi:M\to R_+$ 
fulfilling strict versions of (b02)),
one gets a countable variant of the
Zorn-Bourbaki maximal principle \cite{bourbaki-1949}.
Some conventions are needed.
Let $(<)$ stand 
for the associated relation 
\ben
\item[] 
$x< y$ iff $x\le y$ and $x\ne y$
\een
Clearly, $(<)$ is 
{\it irreflexive} ($x< x$ is false, $\forall x\in M$)
and
{\it transitive} ($x< y$ and $y< z$ imply $x< z$);
it will be referred to as the
{\it strict order} attached to $(\le)$.
Call the point $z\in M$, {\it $(\le)$-maximal}, provided
\ben 
\item[] (b03)\ \ 
$w\in M$, $z\le w$ $\limpl$ $z=w$;\\
or, equivalently:\ 
$M(z,<) (:=\{x\in M; z< x\})$ is empty.
\een

The following (Zorn-Bourbaki) 
maximal version of (BB) 
(denoted, for simplicity, as (BB-Z)) is now available.

\bprop \label{p2}
Suppose that 
the (partially) ordered structure $(M,\le)$ 
is such that
\ben
\item[] (b04)\ \
$(M,<)$ is sequentially inductive: \\
each $(<)$-ascending sequence in $M$
has an upper bound in $M$ (modulo $(<)$)
\item[] (b05)\ \ 
$(M,<)$ is {\it admissible}: \\
there exists at least one function 
$\vphi:M\to R_+$ with the
$(<)$-decreasing  property 
($x< y \limpl \vphi(x)> \vphi(y)$).
\een
Then, $(\le)$ is a Zorn order, in the sense:
for each $u\in M$ there exists a $(\le)$-maximal
$v\in M$ with $u\le v$.
\eprop

\bproof
There are two steps to be passed.

{\bf Step 1.}
We claim that, under these conditions, 
$(M,\le)$ is sequentially inductive.
In fact, let $(x_n; n\ge 0)$ be a 
$(\le)$-ascending sequence in $M$.
If the alternative below is in force
\ben
\item[] 
there exists $k\ge 0$, such that 
$x_k=x_n$, for all $n> k$,
\een
we are done;
because $y:=x_k$ is an upper bound of
$(x_n; n\ge 0)$.
Suppose that 
the opposite alternative is true:
\ben
\item[]
for each $k\ge 0$, there exists 
$h> k$ with $x_k< x_h$.
\een
In this case, we get a 
$(<)$-ascending sequence of ranks $(i(n); n\ge 0)$,
such that 
the subsequence $(y_n:=x_{i(n)}; n\ge 0)$
is $(<)$-ascending.
By the admitted hypothesis, 
there exists $y\in M$ such that 
$y_n< y$, for all $n$.
This, along with the 
$(\le)$-ascending property of $(x_n; n\ge 0)$, 
gives $x_n< y$, for each $n$;
and the claim follows.

{\bf Step 2.}
As $(M,<)$ is admissible,
there exists at least one function 
$\vphi:M\to R_+$ with the
$(<)$-decreasing  property: 
$x< y \limpl \vphi(x)> \vphi(y)$.
Note that, by the very definition of
our strict order $(<)$, 
we have the (converse) representation formula
\ben
\item[]
$x\le y$ iff either $x< y$ or $x=y$.
\een
As a direct consequence of this,
one gets that
$$
\mbox{
$\vphi$ is $(\le)$-decreasing\
($x\le y \limpl \vphi(x)\ge \vphi(y)$).
}
$$
Putting these together,  
(BB) is applicable to $(M,\le)$ and $\vphi$.
From this principle we are assured that, 
given $u\in M$,
there exists a $(\le,\vphi)$-maximal $v\in M$
with $u\le v$.
Suppose by contradiction that 
$v< w$, for some $w\in M$.
As $\vphi$ is $(<)$-decreasing, 
this gives $\vphi(v)> \vphi(w)$;
in contradiction with the
$(\le,\vphi)$-maximal property of $v$.
Hence, $v$ is $(\le)$-maximal;
and we are done.
\eproof

Note that, for the moment,
(BB) $\limpl$ (BB-Z) in 
the {\it strongly reduced axiomatic system} (ZF-AC).
On the other hand, 
this statement includes (see below)
Ekeland's Variational Principle \cite{ekeland-1979} (in short: EVP). 
As a consequence, many extensions of (BB) 
were proposed; 
see, for instance, 
Hyers, Isac and Rassias 
\cite[Ch 5]{hyers-isac-rassias-1997}.
For each (countable) variational principle 
(VP) of this type,  
one therefore has (VP) $\limpl$ (BB) $\limpl$ (EVP);
so, we may ask whether
these inclusions are effective.
As we shall see,
the answer to this is negative.

{\bf (B)}
Let $M$ be a nonempty set; and $\calr\incl M\times M$ 
be a (nonempty) {\it relation} over $M$;
for simplicity, we sometimes write $(x,y)\in \calr$ as $x\calr y$.
Note that $\calr$ may be viewed as a mapping between $M$
and $2^M$ (=the class of all subsets in $M$).
In fact, denote for each $x\in M$ 
\ben
\item[]
$M(x,\calr)=\{y\in M; x\calr y\}$
(=the {\it section} of $\calr$ through $x$);
\een
then, the mapping representation of $\calr$ is
$(\calr(x)=M(x,\calr); x\in M)$.

Call the relation $\calr$ over $M$,
{\it proper} when
\ben
\item[] (b06)\ \ 
$M(c,\calr)$ is nonempty, for each $c\in M$.
\een
Clearly, $\calr$ may be then viewed as a mapping between $M$ 
and $(2)^M$ (=the class of all nonempty subsets in $M$).  

The following "Principle of Dependent Choices" (in short: DC) 
is in effect for our future developments.

\bprop \label{p3}
Suppose that 
$\calr$ is a proper relation over $M$.
Then, for each $a\in M$ there exists 
a sequence $(x_n; n\ge 0)$ in $M$ 
with $x_0=a$ and $x_n\calr x_{n+1}$, for all $n$.
\eprop

This principle, due to 
Bernays \cite{bernays-1942}
and
Tarski \cite{tarski-1948},
is deductible from AC (= the Axiom of Choice), but not conversely;
cf. 
Wolk \cite{wolk-1983}.
Moreover, the 
{\it reduced axiomatic system} (ZF-AC+DC) 
seems to be comprehensive enough for a large part of the
"usual" mathematics; see
Moore \cite[Appendix 2, Table 4]{moore-1982}.
\sk

As an illustration of this assertion,
we show that, ultimately,
(BB) is contained in the underlying
reduced system.

\bprop \label{p4}
We have (DC) $\limpl$ (BB) in 
the strongly reduced system (ZF-AC);
hence, (BB) 
is deductible in 
the reduced system
(ZF-AC+DC).
\eprop

\bproof 
Let the premises of (BB) be admitted;
i.e.:
the quasi-ordered structure
$(M,\le)$ is sequentially inductive
and the function $\vphi:M\to R_+$ 
is $(\le)$-decreasing.
Define the function $\be:M\to R_+$ as:
\ben
\item[]
$\be(v):=\inf[\vphi(M(v,\le))]$,\ $v\in M$.
\een
Clearly, $\be$ is increasing; 
and
\beq \label{201}
\mbox{
$\vphi(v)\ge \be(v)$, for all $v\in M$.
}
\eeq
Moreover, ($\vphi$=decreasing) 
yields a characterization of maximal elements like
\beq \label{202}
\mbox{
$v$ is $(\le,\vphi)$-maximal iff $\vphi(v)=\be(v)$.
}
\eeq
Now, assume by contradiction that the conclusion in this statement 
is false; i.e. [in combination with (\ref{201})+(\ref{202})],
there must be some $u\in M$ such that:
\ben
\item[] (b07)\ \ 
for each $v\in M_u:=M(u,\le)$, one has $\vphi(v)> \be(v)$.
\een
Consequently (for all such $v$), 
$\vphi(v)> (1/2)(\vphi(v)+\be(v)) > \be(v)$;
hence 
\beq \label{203}
\mbox{
$v\le w$ and $(1/2)(\vphi(v)+\be(v))> \vphi(w)$, 
}
\eeq
for at least one $w$ (belonging to $M_u$).
The relation $\calr$ over $M_u$ introduced via (\ref{203})
is proper on $M_u$; i.e.:
\ben
\item[]
$M_u(v,\calr)\ne \es$, for all $v\in M_u$.
\een
So, by (DC), there must be 
a sequence $(u_n)$ in $M_u$ with $u_0=u$ and
\beq \label{204}
\mbox{
$u_n\le u_{n+1}$, $(1/2)(\vphi(u_n)+\be(u_n))> \vphi(u_{n+1})$, 
for all $n$.
}
\eeq
We have thus constructed an ascending sequence  $(u_n)$ in $M_u$
for which the positive (real) sequence $(\vphi(u_n))$ is (via (b07)) 
strictly descending and bounded below; hence 
$\lb:=\lim_n \vphi(u_n)$ exists in $R_+$.
As $(M,\le)$ is sequentially inductive, 
$(u_n)$ is bounded from above in $M$: 
there exists  $v\in M$ such that $u_n\le v$, for all $n$
(whence, $v\in M_u$).
Moreover, since ($\vphi$=decreasing), we must have 
(by the properties of $\be$)

{\bf j)}\
$\vphi(u_n)\ge \vphi(v)$, $\forall n$;\ \
{\bf jj)}\ 
$\vphi(v)\ge \be(v)\ge \be(u_n)$, $\forall n$.

\n
The former of these relations gives $\lb\ge \vphi(v)$ 
(passing to limit as $n\to \oo$).
On the other hand, the latter of these relations yields
(via (\ref{204}))
$$
\mbox{
$(1/2)(\vphi(u_n)+\be(v))> \vphi(u_{n+1})$, for all $n$. 
}
$$
Passing to limit as $n\to \oo$  gives
$$
(\vphi(v)\ge )\be(v)\ge \lb;
$$ 
so, combining with the preceding one,
$$
\mbox{
$\vphi(v)=\be(v)(=\lb)$, contradiction.
}
$$
Hence, (b07) cannot be accepted; and the conclusion follows.
\eproof

Note that, a  slightly different proof
of this may be found in
the 2007  monograph by
C\^{a}rj\u{a} et al 
\cite[Ch 2, Sect 2.1]{carja-necula-vrabie-2007}.
Further metrical aspects of it may be found in 
Turinici \cite{turinici-2002-AUOC}.

{\bf (C)}
In the following, the relationships between (BB)
and 
Ekeland's variational principle \cite{ekeland-1979}
(in short: EVP)
are discussed.

Let $(M,d)$ be a metric space;
and $\vphi:M\to R_+$ be a function.
Assume that 
\ben
\item[] (b08)\ \ 
$(M,d)$ is {\it complete}\  
(each $d$-Cauchy sequence in $M$ is $d$-convergent)
\item[] (b09)\ \ 
$\vphi$ is 
{\it $d$-lsc}:\ 
$\liminf_n \vphi(x_n)\ge \vphi(x)$,\ whenever $x_n\dtends x$; \\
or, equivalently:\
$\{x\in M; \vphi(x)\le t\}$ is $d$-closed, for each $t\in R$.
\een

\bprop \label{p5}
Let these conditions hold. Then,  for each 
(starting point) $u\in M$ 
there exists (another point) $v\in M$ with 
\beq \label{205}
\mbox{
$d(u,v)\le \vphi(u)-\vphi(v)$ (hence $\vphi(u)\ge \vphi(v)$)
}
\eeq
\beq \label{206}
\mbox{
$d(v,x)> \vphi(v)-\vphi(x)$, for each $x\in M\sm \{v\}$.
}
\eeq
\eprop

\bproof 
Let $(\preceq)$ stand for the relation (over $M$):
\ben
\item[]
$x\preceq y$ iff $d(x,y)\le \vphi(x)-\vphi(y)$.
\een
Clearly, $(\preceq)$ acts as a (partial) order
on $M$; note that, as a consequence of this, its
associated relation 
\ben
\item[]
$x\prec y$ iff\ $0< d(x,y)\le \vphi(x)-\vphi(y)$
\een
is a strict order on $X$.
We claim that conditions of (BB-Z) are fulfilled on
$(M,\preceq)$. In fact,  by this very definition, 
$\vphi$ is $(\prec)$-decreasing on $M$;
so that, $(M,\prec)$ is admissible.
On the other hand, let $(x_n)$  be a 
$(\prec)$-ascending sequence in $M$:
\ben
\item[] (b10)\ \ 
$0< d(x_n,x_m)\le \vphi(x_n)-\vphi(x_m)$, if $n< m$.
\een
The sequence $(\vphi(x_n))$ is strictly descending and 
bounded from below; hence a Cauchy one.
This, along with our working hypothesis, 
tells us that $(x_n)$ 
is a $d$-Cauchy sequence in $M$; 
wherefrom by completeness, 
\ben
\item[]
$x_n\dtends y$ as $n\to \oo$, for some $y\in M$.
\een
Passing to limit as $m\to \oo$ in 
the same working hypothesis, one derives 
$$  \mbox{
$d(x_n,y)\le \vphi(x_n)-\vphi(y)$,\
(i.e.: $x_n\preceq y$), for all $n$.
}
$$
This, combined with $(x_n; n\ge 0)$ being $(\prec)$-ascending,
gives $x_n\prec y$, for all $n$;
and shows that 
$(M,\prec)$ is sequentially inductive.
From (BB-Z) it then follows that, for the starting 
$u\in M$ there exists some $v\in M$ with

{\bf h)} $u\preceq v$;\ \   
{\bf hh)} $v\preceq x\in M$ implies $v=x$.

\n
The former of these is just (\ref{205}); 
and the latter one gives at once
(\ref{206}). 
\eproof

This principle found some basic applications to control and optimization,
generalized differential calculus, critical point theory and global
analysis; we refer to the quoted paper  
for a survey of these.
So, it cannot be surprising that, soon after its formulation,
many extensions of (EVP) were proposed.
For example, the dimensional way of extension refers to the
ambient {\it positive halfline} $R_+$ of $\vphi(M)$ being substituted by a
{\it convex) cone} of a
(topological or not) {\it vector space}.
An account of the results in this area is to be found in
the 2003 monograph by 
Goepfert, Riahi, Tammer and Z\u{a}linescu 
\cite[Ch 3]{goepfert-riahi-tammer-zalinescu-2003};
see also 
Turinici \cite{turinici-2002-AUOC}.
On the other hand, the 
{\it (pseudo) metrical} one consists in the conditions
imposed to the ambient metric over $M$ being relaxed.
Some basic results in this direction were obtained by
Kang and Park \cite{kang-park-1990};
see also
Tataru \cite{tataru-1992}. 

{\bf (D)}
By the developments above, 
we therefore have the implications:
$$
\mbox{
(DC) $\limpl$ (BB) $\limpl$ (BB-Z) $\limpl$ (EVP).
}
$$
So, we may ask whether these may be reversed.
Clearly, the natural setting for solving this problem 
is (ZF-AC); referred to (see above) 
as the {\it strongly reduced} Zermelo-Fraenkel system.
\sk

Let $X$ be a nonempty set; 
and $(\le)$ be a (partial) order on it. 
We say that $(\le)$
has  the {\it inf-lattice} property, provided: 
\ben
\item[]
$x\wedge y:=\inf(x,y)$ exists, for all $x,y\in X$.
\een
Remember that $z\in X$ is a $(\le)$-{\it maximal} element if 
$X(z,\le)=\{z\}$; the class of all these points will be
denoted as $\max(X,\le)$. 
Call $(\le)$, a {\it Zorn order} when 
\ben
\item[]
$\max(X,\le)$ is nonempty and {\it cofinal} in $X$ \\
(for each $u\in X$ there exists a $(\le)$-maximal 
$v\in X$ with $u\le v$).
\een
Further aspects are to be described in a metric setting. 
Let $d:X\times X\to R_+$ be a metric over $X$;
and $\vphi:X\to R_+$ be some function.
Then, the natural choice for $(\le)$ above is
\ben
\item[]
$x\le_{(d,\vphi)} y$ iff $d(x,y)\le \vphi(x)-\vphi(y)$;
\een
referred to as the
{\it Br{\o}ndsted order} \cite{brondsted-1976}
attached to $(d,\vphi)$. 
Denote
\ben
\item[]
$X(x,\rho)=\{u\in X; d(x,u)< \rho\}$, $x\in X$, $\rho> 0$
\een
[the {\it open} sphere with center $x$ and radius $\rho$].
Call the ambient metric space  $(X,d)$, {\it discrete} when 
\ben
\item[]
for each $x\in X$ there exists $\rho=\rho(x)> 0$ such that 
$X(x,\rho)=\{x\}$. 
\een
Note that, under such an assumption, 
any function $\psi:X\to R$ is continuous over $X$.
However, the {\it (global) $d$-Lipschitz property} of the same 
\ben
\item[]
$|\psi(x)-\psi(y)|\le L d(x,y)$, $x,y\in X$, for some $L> 0$
\een
cannot be assured, in general.
\sk

Now, the statement below is a particular case of (EVP):

\bprop \label{p6}
Let  the metric space $(X,d)$ and the function $\vphi:X\to R_+$
satisfy
\ben
\item[] (b11)\ \ 
$(X,d)$ is discrete bounded and complete
\item[] (b12)\ \ 
$(\le_{(d,\vphi)})$ has the inf-lattice property
\item[] (b13)\ \ 
$\vphi$ is $d$-nonexpansive and $\vphi(X)$ is countable.
\een
Then, $(\le_{(d,\vphi)})$ is a Zorn order.
\eprop

We shall refer to it as: 
the discrete Lipschitz countable version of EVP 
(in short: (EVP-dLc)).
Clearly, (EVP) $\limpl$ (EVP-dLc). 
The remarkable fact to be added is that
this last principle yields (DC); so, it 
completes the circle between all these.

\bprop \label{p7}
The inclusion below is holding  
(in the strongly reduced Zermelo-Fraenkel system):
(EVP-dLc) $\limpl$ (DC).
So (by the above), 

{\bf i)}
the maximal/variational principles
(BB), (BB-Z) and (EVP) are all equivalent with (DC);
hence, mutually equivalent

{\bf ii)}
each "intermediary" maximal/variational statement (VP)
with (DC) $\limpl$ (VP) $\limpl$ (EVP)
is equivalent with both (DC) and (EVP).
\eprop

For a complete proof, see
Turinici \cite{turinici-2011-AUAICI}.
In particular, when the 
discrete, bounded, inf-lattice and nonexpansive 
properties are ignored in (EVP-dLc), 
the last result above reduces to the one in 
Brunner \cite{brunner-1987}.
Note that, in the same particular setting,
a different proof of (EVP) $\limpl$ (DC)
was provided in 
Dodu and Morillon \cite{dodu-morillon-1999}.
Further aspects may be found in 
Schechter \cite[Ch 19, Sect 19.51]{schechter-1997}.

\section{Cantor complete ultrametrics}
\setcounter{equation}{0}

Let $X$ be a nonempty set.
By an {\it ultrametric}
(or: {\it non-Archimedean metric}) 
on $X$, we mean any mapping
$d:X\times X\to R_+$ 
with the properties:
\ben
\item[] (c01)\ \ 
$x=y$ iff $d(x,y)=0$ \ \ 
\hfill (reflexive sufficient)
\item[] (c02)\ \ 
$d(x,y=d(y,x)$, $\forall x,y\in X$\ \
\hfill (symmetric)
\item[] (c03)\ \ 
$d(x,z)\le \max\{d(x,y),d(y,z)\}$,
$\forall x,y,z\in X$ \ \
\hfill (ultra-triangular);
\een
in this case, the pair $(X,d)$ will be 
referred to as an {\it ultrametric space}.
Note that, any 
ultrametric is a
(standard) metric (on $X$), because
$$
d(x,z)\le \max\{d(x,y),d(y,z)\}\le
d(x,y)+d(y,z),\ 
\forall x,y,z\in X;
$$
but, the converse 
is not in general valid.
The class of 
these ultrametrics is nonempty.
In fact, the {\it discrete} metric on $X$
introduced as: for each $x,y\in X$
\ben
\item[]
$d(x,y)=1$ if $x\ne y$;\
$d(x,y)=0$, if $x=y$,
\een
is an ultrametric, as it can be directly
seen.
Further examples may be found in 
Rooij \cite[Ch 3]{rooij-1978}.

Let in the following 
$(X,d)$ be an ultrametric space.
Note that,
the presence of ultra-triangular inequality
induces a lot of dramatic changes 
with respect to the standard metrical case;
some basic ones will 
be shown below.
[These were stated without proof in
Khamsi and Kirk
\cite[Ch 5, Sect 5.7]{khamsi-kirk-2001};
see also 
Rooij \cite[Ch 2]{rooij-1978};
however, for completeness 
reasons, we shall provide a proof
of them].

\blemma \label{le1}
Let $x,y,z\in X$ be such that 
$d(x,y)\ne d(y,z)$.
Then, necessarily,
$$
d(x,z)=\max\{d(x,y),d(y,z)\};
$$
hence, either $d(x,z)=d(x,y)$
or $d(x,z)=d(y,z)$.
In other words:
each triangle $(x,y,z)$ in $X$
is $d$-isosceles.
\elemma

\bproof
Suppose by contradiction that 
\ben
\item[]
$d(x,z)< \max\{d(x,y),d(y,z)\}$.
\een
We have two alternatives to consider:

{\bf i)}
Suppose that $d(x,y)<d(y,z)$.
By the working hypothesis,
we then have 
$d(x,z)< d(y,z)$.
In this case, the ultra-triangular 
inequality gives
$$
d(y,z)\le \max\{d(x,y),d(x,z)\}< d(y,z);\
\mbox{contradiction}.
$$

{\bf ii)}
Suppose that 
$d(y,z)< d(x,y)$.
By the working hypothesis,
we then have
$d(x,z)< d(x,y)$;
so, again from the
ultra-triangular inequality,
$$
d(x,y)\le \max\{d(x,z),d(y,z)\}< d(x,y);
\mbox{contradiction}.
$$

Having discussed all possible 
alternatives, we are done.
\eproof
  
By definition, any set of the form
\ben
\item[] (c04)\ \ 
$X[a,r]=\{x\in X; d(a,x)\le r\}$,\
$a\in X$, $r\in R_+$,
\een
will be referred to as a 
{\it $d$-closed sphere} with center $a\in X$
and radius $r\in R_+$;
note that this is a nonempty subset of $X$,
in view of $a\in X[a,r]$. 
In the following, some results involving the
family of all $d$-closed spheres
\ben
\item[] 
$\calm=\{X[a,r]; a\in X, r\in R_+\}\incl (2)^X$
\een
will be discussed.

\blemma \label{le2}
Let $M_1:=X[a_1,r_1]$, $M_2:=X[a_2,r_2]$
be a couple of non-disjoint $d$-closed spheres in $X$. 
Then,
\ben
\item[] (i)\ \ 
$M_1\incl M_2$, whenever $r_1\le r_2$
\item[] (ii)\ \ 
$M_1=M_2$, whenever $r_1=r_2$.
\een
\elemma

\bproof
As $M_1\cap M_2\ne \es$,
there exists at least one element 
$b\in M_1\cap M_2$. 

i):\
Assume that $r_1\le r_2$;
and let $x\in M_1$ be arbitrary fixed.
From the ultra-triangular inequality, we have
$$
d(x,b)\le \max\{d(x,a_1),d(b,a_1)\}\le r_1;
$$
and this in turn yields
(by the same procedure)
$$
d(x,a_2)\le \max\{d(x,b),d(b,a_2)\}\le 
\max\{r_1,r_2\}=r_2;\
\mbox{i.e.:}\ x\in M_2.
$$
As $x\in M_1$ was arbitrarily chosen,
one derives $M_1\incl M_2$.

ii):\
Evident, by the preceding step.
\eproof

\blemma \label{le3}
Let $a,b\in X$ and $s\ge 0$
be such that $a\in X[b,s]$.
Then,
\beq \label{301}
\mbox{
$X[a,r]\incl X[b,s]$, for each 
$r\in[0,s]$.
}
\eeq
\elemma

\bproof
Let $r\in [0,s]$ be arbitrary fixed.
By the imposed hypothesis,
$$
\mbox{
$a\in X[a,r]\cap X[b,s]$;\
whence, 
$X[a,r]\cap X[b,s]\ne \es$;
}
$$
and then, from the previous result,
we are done.
\eproof

The next statement is, in a certain sense,
a reciprocal of the previous one.
Denote 
\ben
\item[]
($Y_1,Y_2\in 2^X$):\ 
$Y_1\sincl Y_2$ iff $Y_1\incl Y_2$
and $Y_1\ne Y_2$.
\een
Clearly, $(\sincl)$ is nothing else 
than the 
{\it strict order} 
(i.e.: irreflexive and transitive relation) 
attached to the usual (partial) order 
$(\incl)$ over $2^X$.

\blemma \label{le4}
Let $M_1:=X[a_1,r_1]$, $M_2:=X[a_2,r_2]$
be two $d$-closed balls in $X$. Then,
\beq \label{302}
M_1\sincl M_2 \limpl r_1< r_2.
\eeq
\elemma

\bproof
Suppose that $M_1\sincl M_2$;
but (contrary to the conclusion) $r_2\le r_1$.
As $M_1\cap M_2=M_1\ne \es$,
one has by a preceding result
(and the working hypothesis)
$$
\mbox{
$M_2\incl M_1\sincl M_2$;\ contradiction.
}
$$
This proves our assertion.
\eproof

We are now introducing a basic notion.
Call the ultrametric space
$(X,d)$, 
{\it Cantor strongly complete}
(in short: {\it Cantor s-complete}),
provided
\ben
\item[] (c05)\ \ 
each $(\aincl)$-ascending sequence
$(M_n:=X[a_n,r_n]; n\ge 0)$ in $\calm$ \\
has a nonempty intersection.
\een
A (formally) weaker variant of this definition is
as follows.
Call the ultrametric space
$(X,d)$, 
{\it Cantor complete}, 
provided
\ben
\item[] (c06)\ \ 
each $(\saincl)$-ascending sequence
$(M_n:=X[a_n,r_n]; n\ge 0)$ in $\calm$ \\
has a nonempty intersection.
\een
Clearly, we have
\beq \label{303} 
\mbox{
($\forall$ ultrametric structure):
Cantor s-complete $\limpl$ Cantor complete.
}
\eeq
The reciprocal inclusion is also true,
as results from

\blemma \label{le5}
For each ultrametric structure $(X,d)$, 
we have
\beq \label{304}
\mbox{
\btab{l}
Cantor complete $\limpl$ Cantor s-complete; \\
hence,
Cantor complete $\lequi$ Cantor s-complete.
\etab
}
\eeq
\elemma

\bproof
Suppose that the ultrametric space 
$(X,d)$ is Cantor complete;
and let $(M_n:=X[a_n,r_n]; n\ge 0)$
be a $(\aincl)$-ascending sequence in 
$\calm$.
If one has that
\ben
\item[]
$\exists (i\ge 0)$, $\forall (j> i)$:\ $M_i=M_j$,
\een
we are done; because
$\cap\{M_n; n\ge 0\}=M_i$.
Suppose now that the opposite 
alternative is holding:
\ben
\item[]
$\forall (i\ge 0)$, $\exists (j> i)$:\ $M_i\saincl M_j$.
\een
There exists then a strictly ascending sequence of
ranks $(i(n); n\ge 0)$, such that 
the subsequence 
$(L_n:=M_{i(n)}; n\ge 0)$ of $(M_n; n\ge 0)$ fulfills
$$
\mbox{
$(L_n)$ is $(\saincl)$-ascending:\
$p< q$ $\limpl$ $L_p \saincl L_q$.
}
$$
By the imposed hypothesis,
$L:=\cap\{L_n; n\ge 0\}$ is nonempty.
This, along with
$L=\cap\{M_n; n\ge 0\}$, ends the
argument.
\eproof

Denote, for simplicity
\ben
\item[]
$\Ga=X\times R_+$;\ 
hence, $\Ga=\{(a,\rho); a\in X, \rho\in R_+\}$.
\een
A natural relation to be
introduced here is the following
\ben
\item[]
$(a,\rho)\prec (b,\sig)$ iff  $X[a,\rho]\saincl X[b,\sig]$.
\een
Clearly, $(\prec)$ is
irreflexive and transitive;
hence, a strict order on $\Ga$.
Let $(\preceq)$ stand for the associated
(partial) order
\ben
\item[]
$(a,\rho)\preceq (b,\sig)$ iff 
either $(a,\rho)\prec (b,\sig)$ 
or 
$(a,\rho)=(b,\sig)$.
\een
Having these precise, let us 
introduce the function 
\ben
\item[]
($\vphi:\Ga\to R_+$):\
$\vphi(a,\rho)=\rho$,\ $(a,\rho)\in \Ga$.
\een
By a previous result, we have
\ben
\item[]
$\vphi$ is $(\prec)$-decreasing:\
$(a,\rho)\prec (b,\sig)$ 
$\limpl$ $\vphi(a,\rho)> \vphi(b,\sig)$.
\een
This tells us that, necessarily,
\beq \label{305}
\mbox{
$(\Ga,\prec)$ is admissible;\
hence, so is $(\De,\prec)$, 
where $\es\ne \De\incl \Ga$.
}
\eeq
As a consequence, the following 
(relative) maximal result 
is available.

\bprop \label{p8}
Let the (nonempty) subset $\De$ of $\Ga$
be such that
\ben
\item[]
$(\De,\prec)$ is sequentially inductive:\\
each $(\prec)$-ascending sequence in $\De$ is bounded 
above in $\De$ (modulo $(\prec)$).
\een
Then, $(\preceq)$ is a Zorn order on $\De$; 
i.e.:
for each (starting element) 
$(a,\rho)\in \De$, there exists 
(another element) $(b,\sig)\in \De$, with 
\ben
\item[] i)\ \ 
$(a,\rho)\preceq (b,\sig)$; i.e.:
either $(a,\rho)\prec (b,\sig)$ 
or
$(a,\rho)=(b,\sig)$
\item[] ii)\ \ 
$(b,\sig)\prec (c,\tau)$ is impossible,
for each $(c,\tau)\in \De$.
\een
\eprop

\bproof
By the admissible property for $\Ga$,
we have 
$$
\mbox{
(the strictly ordered structure) 
$(\De,\prec)$ is admissible.
}
$$
Combining with the admitted hypothesis,
it results that 
the sequential type 
maximal result (BB-Z) 
is applicable to $(\De,\preceq)$;
and, from this, we are done.
\eproof

\section{Application (fixed point theorems)}
\setcounter{equation}{0}

In the following, an application of the above
developments is given
to the ultrametric fixed point theory.

Let $(X,d)$ be an ultrametric space.
We say that $T\in \calf(X)$ is 
{\it $d$-strictly-nonexpansive}, provided
\ben
\item[] (d01)\ \ 
$d(Tx,Ty)< d(x,y)$,\ $\forall x,y\in X$,\ $x\ne y$.
\een
Note that, in particular, $T$ is {\it $d$-nonexpansive}:
\ben
\item[] (d02)\ \  
$d(Tx,Ty)\le d(x,y)$,\ for all $x,y\in X$.
\een
The following fixed point theorem 
over ultrametric spaces is available.

\btheorem \label{t2}
Suppose that 
$T$ is $d$-strictly-nonexpansive (see above).
In addition, let $(X,d)$ be Cantor complete. 
Then, $T$ is fix-singleton; 
whence, it has a unique fixed point in $X$.
\etheorem

\bproof
There are several steps to be followed.

{\bf Step 1.}
By the $d$-strict-nonexpansive property,
we have
$$
\mbox{
$\Fix(T)$ is asingleton; 
i.e.: $T$ is fix-asingleton.
}
$$
So, all we have to establish is that 
$\Fix(T)$ appears as nonempty.

{\bf Step 2.}
Remember that, over $\Ga:=X\times R_+$
we introduced the strict ordering
\ben
\item[]
$(a,\rho)\prec (b,\sig)$ iff  $X[a,\rho]\saincl X[b,\sig]$;
\een
as well as the associated ordering 
\ben
\item[]
$(a,\rho)\preceq (b,\sig)$ iff 
either $(a,\rho)\prec (b,\sig)$ 
or 
$(a,\rho)=(b,\sig)$.
\een
Moreover, we have that 
\beq \label{401}
\mbox{
$(\Ga,\prec)$ is admissible;
hence, so is $(\De,\prec)$, where
$\es\ne \De\incl \Ga$.
}
\eeq

{\bf Step 3.}
Denote, for simplicity
\ben
\item[]
$\De=\{(a,d(a,Ta)); a\in X\}$;
\een
this is a nonempty subset of $\Ga$.
By a previous relation, we have that 
\beq \label{402}
\mbox{
(the strictly ordered structure) 
$(\De,\prec)$ is admissible.
}
\eeq
Moreover, we claim that 
the structure 
$(\De,\prec)$ is sequentially inductive.
In fact, let 
$((a_n,d(a_n,Ta_n)); n\ge 0)$
be a $(\prec)$-ascending 
sequence in $\De$; i.e.:
\ben
\item[]
$M_i\saincl M_j$,\ for $i< j$;\
where
$(M_n:=X[a_n,d(a_n,Ta_n)]; n\ge 0)$.
\een
As $(X,d)$ is Cantor complete, it follows that
$$
\mbox{
$L:=\cap\{M_n; n\ge 0\}$ is nonempty;
}
$$
let $b\in L$ be some point of it.
By the very definition above
(and the $d$-nonexpansive property of $T$)
$$
d(Tb,Ta_n)\le d(b,a_n)\le d(a_n,Ta_n),\ 
\forall n\ge 0.
$$
Combining with the ultra-triangular inequality,
one gets (for the same ranks)
$$
d(b,Tb)\le \max\{d(b,a_n),d(a_n,Ta_n),d(Ta_n,Tb)\}
\le d(a_n,Ta_n);
$$
and this, by a previous auxiliary fact, yields
$$
\mbox{
$X[a_n,d(a_n,Ta_n)]\saincl X[b,d(b,Tb)]$,\ for all $n$;
}
$$
or, equivalently (by definition)
$$
\mbox{
$(a_n,d(a_n,Ta_n))\prec (b,d(b,Tb))\in \De$, for all $n$; 
}
$$
which proves the desired fact.

{\bf Step 4.}
Putting these together, it follows that
the previous maximal principle
is applicable to $(\De,\preceq)$.
So, for the starting element
$(u,d(u,Tu))$ in $\De$,
there exists another element
$(v,d(v,Tv))$ in $\De$, with
\ben
\item[] i)\ \ 
$(u,d(u,Tu))\preceq (v,d(v,Tv))$
\item[] ii)\ \ 
for each $w\in X$, 
$(v,d(v,Tv))\prec (w,d(w,Tw))$ is impossible.
\een
Suppose by contradiction that
\ben
\item[]
$d(v,Tv)> 0$;
hence, 
$d(Tv,T^2v)< d(v,Tv)$.
\een
We claim that
\ben
\item[]
$(v,d(v,Tv))\prec (Tv,d(Tv,T^2v))$;
\een
and this, by the previous 
maximal property of $(v,d(v,Tv))$
yields a contradiction.
The desired relation may be written as
\ben
\item[]
$X[v,d(v,Tv)]\saincl X[Tv,d(Tv,T^2v)]$;
\een
to establish it, we may 
proceed as follows.

{\bf I)}
Let $y\in X[Tv,d(Tv,T^2v)]$ be arbitrary fixed;
hence, 
\ben
\item[]
$d(y,Tv)\le d(Tv,T^2v) (< d(v,Tv))$.
\een
By the ultra-triangular inequality,
$$
d(y,v)\le \max\{d(y,Tv),d(v,Tv)\}
= d(v,Tv); 
$$
whence, $y\in X[v,d(v,Tv)]$;
this, by the arbitrariness of $y$, gives 
$$
X[Tv,d(Tv,T^2v)]\incl X[v,d(v,Tv)].
$$

{\bf II)}
From the working assumption about $v$, 
one must have
$$
\mbox{
($v\in X[v,d(v,Tv)]$ and)\ $v\notin X[Tv,d(Tv,T^2v)]$;
}
$$
hence, the above inclusion is strict.
The proof is thereby complete.
\eproof

By the argument above, 
this fixed point result is
a consequence of the Brezis-Browder 
ordering principle
\cite{brezis-browder-1976};
hence, ultimately, it 
is deductible in the reduced Zermelo-Fraenkel
system (ZF-AC+DC).
Note that, similar conclusions
are to be derived  
for the related fixed point results 
over ultrametric spaces
due to 
Gaji\'{c} \cite{gajic-2001}
and 
Pant \cite{pant-2014};
see also 
Wang and Song \cite{wang-song-2013}.
Further aspects of this theory
concerning fuzzy ultrametric spaces
may be found in
Sayed \cite{sayed-2013}.



\begin{thebibliography}{99}



\bibitem{banach-1922}
{S. Banach}, 
\it Sur les op\'{e}rations dans les ensembles abstraits 
et leur application aux \'{e}quations int\'{e}grales,
\rm Fund. Math., 3 (1922), 133-181.


\bibitem{bernays-1942}
{P. Bernays},
\it A system of axiomatic set theory: 
Part III. Infinity and enumerability analysis,
\rm J. Symbolic Logic, 7 (1942), 65-89.


\bibitem{bourbaki-1949}
{N. Bourbaki},
\it Sur le th\'{e}or\`{e}me de Zorn,
\rm Archiv Math., 2 (1949/1950), 434-437.


\bibitem{boyd-wong-1969}
{D. W. Boyd} and {J. S. W. Wong},
\it On nonlinear contractions, 
\rm Proc. Amer. Math. Soc., 20 (1969), 458-464. 


\bibitem{brezis-browder-1976}
{H. Brezis} and {F. E. Browder},
\it A general principle on ordered sets in nonlinear functional analysis,
\rm Advances Math., 21 (1976), 355-364.


\bibitem{brondsted-1976}
{A. Br{\o}ndsted},
\it Fixed points and partial orders,
\rm Proc. Amer. Math. Soc., 60 (1976), 365-366.


\bibitem{brunner-1987}
{N. Brunner},
\it Topologische Maximalprinzipien,
\rm Zeitschr. Math. Logik Grundl. Math., 33 (1987), 135-139.


\bibitem{carja-necula-vrabie-2007}
{O. C\^{a}rj\u{a}}, {M. Necula} and {I. I. Vrabie},
\it Viability, Invariance and Applications,
\rm North Holland Math. Studies vol. 207, Elsevier B. V., Amsterdam, 2007.


\bibitem{cohen-1966}
{P. J. Cohen},
\it Set Theory and the Continuum Hypothesis, 
\rm Benjamin, New York, 1966.


\bibitem{dodu-morillon-1999}
{J. Dodu} and {M. Morillon},
\it The Hahn-Banach property and the Axiom of Choice,
\rm Math. Logic Quarterly, 45 (1999), 299-314.


\bibitem{edelstein-1961}
{M. Edelstein},
\it An extension of Banach's contraction principle,
\rm Proc. Amer. Math. Soc., 12 (1961), 7-10.


\bibitem{ekeland-1979}
{I. Ekeland},
\it Nonconvex minimization problems,
\rm Bull. Amer. Math. Soc. (New Series), 1 (1979), 443-474.


\bibitem{gajic-2001}
{L. Gaji\'{c}},
\it On ultrametric spaces,
\rm Novi Sad J. Math., 31 (2001), 69-71.


\bibitem{goepfert-riahi-tammer-zalinescu-2003}
{A. Goepfert}, {H. Riahi}, {C. Tammer} and {C. Z\u{a}linescu},
\it Variational Methods in Partially Ordered Spaces,
\rm Canad. Math. Soc. Books Math. vol. 17, Springer, New York, 2003.


\bibitem{hyers-isac-rassias-1997}
{D. H. Hyers}, {G. Isac} and {T. M. Rassias},
\it Topics in Nonlinear Analysis and Applications,
\rm World Sci. Publ., Singapore, 1997.


\bibitem{kang-park-1990}
{B. G. Kang} and {S. Park},
\it On generalized ordering principles in nonlinear analysis,
\rm Nonlin. Anal., 14 (1990), 159-165.


\bibitem{khamsi-kirk-2001}
{M. A. Khamsi} and {W. A. Kirk},
\it An Introduction to Metric Spaces and Fixed Point Theory,
\rm John Wiley \& Sons Inc., New York, 2001.
 

\bibitem{leader-1979}
{S. Leader},
\it Fixed points for general contractions in metric spaces,
\rm Math. Japonica, 24 (1979), 17-24.

 
\bibitem{matkowski-1975}
{J. Matkowski},     
\it Integrable solutions of functional equations, 
\rm Dissertationes Math., Vol. 127, 
Polish Sci. Publ., Warsaw, 1975.


\bibitem{mishra-pant-2014}
{S. N. Mishra} and {R. Pant},
\it  Generalization of some fixed point theorems in 
ultrametric spaces,
\rm Adv. Fixed Point Th., 4 (2014), 41-47.


\bibitem{moore-1982}
{G. H. Moore},
\it Zermelo's Axiom of Choice: its Origin, Development and Influence,
\rm Springer, New York, 1982.


\bibitem{pant-2014}
{R. Pant},
\it Some new fixed point theorems for contractive 
and nonexpansive mappings,
\rm Filomat, 28 (2014), 313-317.


\bibitem{petalas-vidalis-1993}
{C. Petalas} and {T. Vidalis},
\it A fixed point theorem in non-Archimedean vector spaces,
\rm Proc. Amer. Math. Soc., 118 (1993), 819-821.
 

\bibitem{rhoades-1977}
{B. E. Rhoades},
\it A comparison of various definitions of contractive mappings,
\rm Trans. Amer. Math. Soc., 226 (1977), 257-290.


\bibitem{rooij-1978}
{A. C. M. van Rooij},
\it Non-Archimedean Functional Analysis,
\rm Marcel Dekker Inc., New York, 1978.


\bibitem{rus-2001}
{I. A. Rus}, 
\it Generalized Contractions and Applications, 
\rm Cluj University Press, Cluj-Napoca, 2001.


\bibitem{sayed-2013}
{A. F. Sayed},
\it Common fixed point theorems of multivalued maps
in fuzzy ultrametric spaces,
\rm J. Math., Volume 2013, Article ID 617532.


\bibitem{schechter-1997}
{E. Schechter},
\it Handbook of Analysis and its Foundation,
\rm Academic Press, New York, 1997. 


\bibitem{tarski-1948}
{A. Tarski},
\it Axiomatic and algebraic aspects of two theorems on sums of cardinals,
\rm Fund. Math., 35 (1948), 79-104.


\bibitem{tataru-1992}
{D. Tataru},
\it Viscosity solutions of Hamilton-Jacobi equations with unbounded
nonlinear terms,
\rm J. Math. Anal. Appl., 163 (1992), 345-392.


\bibitem{turinici-1986-DM}
{M. Turinici},
\it Fixed points for monotone iteratively local contractions,
\rm Dem. Math., 19 (1986), 171-180.


\bibitem{turinici-2002-AUOC}
{M. Turinici},
\it Minimal points in product spaces,
\rm An. \c{S}t. Univ. "Ovidius" Constan\c{t}a (Ser. Mat.), 10 (2002), 109-122.


\bibitem{turinici-2011-AUAICI}
{M. Turinici},
\it  Brezis-Browder principle and Dependent Choice,
\rm An. \c{S}t. Univ. "Al. I. Cuza" Ia\c{s}i (S. N.), Mat., 57 (2011), 263-277.


\bibitem{wang-song-2013}
{Q. Wang} and {M. Song},
\it Some coupled fixed point theorems in ultra metric spaces,
\rm Sci. J. Math. Res., 3 (2013), 114-118.


\bibitem{wolk-1983}
{E. S. Wolk},
\it On the principle of dependent choices and some forms of Zorn's lemma,
\rm Canad. Math. Bull., 26 (1983), 365-367.



\end{thebibliography}
\end{document}